\documentclass[a4paper]{article}
\usepackage{amsfonts}
\usepackage{amsmath}
\usepackage{amssymb}
\begin{document}
\newtheorem{theorem}{Theorem}[section]
\newtheorem{lemma}[theorem]{Lemma}
\newtheorem{corollary}[theorem]{Corollary}
\newtheorem{conjecture}[theorem]{Conjecture}
\newtheorem{remark}[theorem]{Remark}
\newtheorem{definition}[theorem]{Definition}
\newtheorem{problem}[theorem]{Problem}
\newtheorem{example}[theorem]{Example}
\newtheorem{proposition}[theorem]{Proposition}
\title{{\bf Variation of Bergman kernels of adjoint line bundles}}
\date{November 13, 2005}
\author{Hajime TSUJI}
\maketitle
\begin{abstract}
Let $f : X \longrightarrow S$ be a smooth projective family and let 
$(L,h)$ be a singular hermitian line bundle on $X$ with semipositive curvature
current.  

Let $K_{s}:= K(X_{s},K_{X_{s}} + L\mid X_{s},h\mid X_{s}) (s\in S)$ 
be the Bergman kernel of $K_{X_{s}} + L\mid X_{s}$ with respect to 
$h\mid X_{s}$ and let  $h_{B}$ the singular hermitian metric 
on $K_{X}+L$ defined by $h_{B}\mid_{X_{s}}:= 1/K_{s}$.  
We prove that $h_{B}$ has semipositive curvature.
This is a generalization of the recent result of Berndtsson (\cite{b}).

Using this result, we give a new proof of Kawamata's semipositivity theorem for 
the direct image of relative multi canonical bundle. 
\end{abstract}
\tableofcontents
\section{Introduction}

The theory of Bergman kernels was initiated by S. Bergman (\cite{berg})
in 1933.
But the variation of Bergman kernels has not  been studied 
until quite recently. 
  In fact in 2004, F. Maitani and H. Yamaguchi proved the following theorem 
and initiated the study of the variation of Bergman kernels. 

\begin{theorem}(\cite{m-y})\label{m-y} Let $\Omega$ be a pseudoconvex domain 
in $\mbox{\bf C}_{z}\times \mbox{C}_{w}$ with smooth boundary. 
Let $\Omega_{t} := \Omega \cap (\mbox{\bf C}_{z}\times \{ t\})$ and 
Let $K(z,t)$ be the Bergman kernel function of 
$\Omega_{t}$. 

Then $\log K(z,t)$ is a plurisubharmonic function on $\Omega$. $\square$
\end{theorem}
Recently  generalizing Theorem \ref{m-y}, B. Berndtsson proved the following
higher dimensional and twisted version of Theorem \ref{m-y}.
 
\begin{theorem}(\cite{b})\label{b}
Let $D$ be a pseudoconvex domain in $\mbox{\bf C}^{n}_{z}\times \mbox{\bf C}^{k}_{t}$.   And let $\phi$ be a plurisubharmonic function on $D$.
For $t\in \Delta$, we set  $D_{t} := \Omega \cap (\mbox{\bf C}^{n}\times \{ t\})$
and $\phi_{t} := \phi\mid D_{t}$. 
Let $K(z,t) (t\in \mbox{\bf C}^{k}_{t})$ be the Bergman kernel of the Hilbert space 
\[
A^{2}(D_{t},e^{-\phi_{t}}) := \{ f\in {\cal O}(\Omega_{t})\mid 
\int_{D_{t}}e^{-\phi_{t}}\mid f\mid^{2} < + \infty \} .
\]
Then $\log K(z,t)$ is a plurisubharmonic function on $D$. $\square$
\end{theorem}
As in mensioned in \cite{b2}, his proof also works for a pseudoconvex 
domain in a locally trivial family of manifolds which admits 
a Zariski dense Stein subdomain.  

Using Theorem \ref{b}, he proved the following theorem. 

\begin{theorem}(\cite[Theorem 1.1]{b2})\label{b2}
Let us consider a domain $D = U \times\Omega$ and let $\phi$ be a 
plurisubharmonic function on $D$.
For simplicity we assume that $\phi$ is smooth up to the boundary and 
strictly plurisubharmonic in $D$. 
Then for each $t\in U$, $\phi_{t} := \phi (\cdot ,t)$ is plurisubharmonic 
on $\Omega$.  Let $A^{2}_{t}$ be the Bergman space of holomorphic functions
on $\Omega$ with norm
\[
\parallel f\parallel^{2} = \parallel f\parallel^{2}_{t}:=  \int_{\Omega} e^{-\phi_{t}}\mid f\mid^{2}. 
\]
The spaces $A^{2}_{t}$ are all equal as vector spaces but have norms 
that vary with $t$. 
Then ``infinite rank'' vector bundle $E$ over $U$ with fiber
$E_{t} = A^{2}_{t}$ is therefore trivial as a bundle but is equipped with 
a notrivial metric. 
Then $(E,\parallel \,\,\,\,\parallel_{t})$ is strictly positive 
in the sense of Nakano. $\square$ 
\end{theorem}
In Theorem \ref{b} the assumption that $D$ is a pseudoconvex domain 
in the product space is rather strong.  And in Theorem \ref{b2}, Berndtsson
also assumed that $D$ is a product. 

In this paper we shall remove these assumptions and generalize Theorems \ref{b},\ref{b2} to the case of adjoint line bundles smooth projective fibrations. 

By using this generalization we can  study nonlocally trivial algebraic 
fiber space. 

To state our theorem, let us introduce the notion of the Bergman kernels of adjoint line bundles. 
Let $X$ be a complex manifold of dimension $n$ and let $(L,h)$ be a singular hermitian line bundle (cf. Definition \ref{singhm}) on $X$. 
Let $K_{X}$ denote the canonical line bundle on $X$. 
Let $A^{2}(X,K_{X}+L,h)$ be the Hilbert space defined by
\[
A^{2}(X,K_{X}+L,h) := \{ \sigma \in H^{0}(X,{\cal O}_{X}(K_{X}+L)) \mid
(\sqrt{-1})^{n^{2}}\int_{X}h \cdot \sigma \wedge \bar{\sigma} < + \infty \} ,
\]
where we have defined the inner product on $A^{2}(X,K_{X}+L,h)$ by 
\[
(\sigma ,\tau ):= (\sqrt{-1})^{n^{2}}\int_{X}h \cdot \sigma \wedge \bar{\tau} .
\]
We define the Bergman kernel $K(X,K_{X}+ L,h)$ of the adjoint bundle $K_{X} + L$ with respect to $h$ by 
\[
K(X,K_{X}+ L,h) := (\sqrt{-1})^{n^{2}}\sum_{i}\sigma_{i}\wedge\bar{\sigma}_{i}.
\]
where $\{ \sigma_{i}\}$ is a complete orthonormal basis of the Hilber space
$A^{2}(X,K_{X}+L,h)$. 
Then $K(X,K_{X}+ L,h)$ is independent of the choice of the complete orthonormal
basis $\{\sigma_{i}\}$. 
In fact 
\[
K(X,K_{X}+ L,h)(x) = \sup \{ (\sqrt{-1})^{n^{2}}\sigma (x)\wedge\bar{\sigma}(x) \mid
\parallel\sigma\parallel = 1\}
\]
holds. \vspace{5mm}\\
Now we shall state the main theorem in this paper. 

\begin{theorem}\label{v}
Let $f : X \longrightarrow S$ be a smooth projective family 
of projective varieties over a complex manifold $S$.
Let $(L,h)$ be a singular hermitian line bundle on $X$ such that 
$\Theta_{h}$ is semipositive on $X$.
Let $K_{s} := K(X_{s},K_{X}+ L\mid_{X_{s}},h\mid_{X_{s}})$ be the Bergman kernel 
of $K_{X_{s}}+ (L\mid X_{s})$ with respect to $h\mid X_{s}$.
Then the singular hermitian metric $h_{B}$ of $K_{X/S} + L$ defined by 
\[
h_{B}\mid X_{s}:= K_{s}^{-1}
\]
has semipositive curvature on $X$. $\square$ 
\end{theorem}
Theorem \ref{v} follows from Theorem \ref{b} by a simple trick as follows. 
We may assume that $S$ is the unit open disk $\Delta$ cetered at $O$. 
$f : X \longrightarrow S$ is not locally trivial. 
We shall embed $X$ into the trivial family $p : X \times \Delta \longrightarrow \Delta
,p(x,t) = x (x\in X, t\in \Delta )$
by 
\[
i : X \longrightarrow X \times \Delta 
\]  
defined by 
\[
i(x) := (x,f(x)). 
\]
Then $i(X)$ is a hypersurface in $X\times \Delta$ and not a domain in $X\times\Delta$. 
So we shall thicken $i(X)$ by replacing $X_{t} (t\in \Delta )$ 
by $f^{-1}(\Delta (t,\varepsilon ))$, where $\Delta (t,\varepsilon )$ 
denotes the open disk of radius $\varepsilon$ centered at $t$. 
In this way we construct a thickend family 
\[
f_{\varepsilon} : X(\varepsilon )\longrightarrow \Delta (1/2)
\]
which is considered to be a pseudoconvex domain in the product family $X\times\Delta (1/2)$ over $\Delta (1/2)$,
where $\Delta (1/2)$ denotes $\Delta (0,1/2)$. 
Then Theorem \ref{b} is applicable to the family of Bergman kernels 
of the adjoint bundle of $p^{*}(L,h)$ over $\Delta (1/2)$. 
Letting $\varepsilon$ tend to $0$, with the rescaling constant $\pi\varepsilon^{2}$, we obtain Theorem \ref{v}. 

As a direct consequence, we can also generalize Theorem \ref{b2} as follows. 
 
\begin{theorem}\label{nakano}
Let $f : X \longrightarrow S$ be a smooth projective family 
of over a complex curve $S$ of relative dimension $n$.
Let $(L,h)$ be a hermitian line bundle on $X$ such that 
$\Theta_{h}$ is semipositive on $X$.
We define the hermitian metric $h_{E}$ on  $E:= f_{*}{\cal O}_{X}(K_{X/S} +L)$
 by 
\[
h_{E}(\sigma ,\tau ) := (\sqrt{-1})^{n^{2}}\int_{X_{s}}h\cdot \sigma\wedge\bar{\tau}.
\] 
Then $(E,h_{E})$ is semipositive in the  sense of Nakano. 
Moreover if $\Theta_{h}$ is strictly positive, then 
$(E,h_{E})$ is strictly positive in the sense of Nakano. $\square$ 
\end{theorem}

\begin{theorem}\label{kawamata}(\cite[p.57, Theorem 1]{ka1})
Let $f : X \longrightarrow C$ be an algebraic fiber space over a projective 
curve $C$. 
Then $f_{*}{\cal O}_{X}(mK_{X/S})$ is a semipositive vector bundle on $S$
in the sense that for  any quotient sheaf ${\cal Q}$ of $f_{*}{\cal O}_{X}(mK_{X/S})$, $\deg_{C}{\cal Q} \geqq 0$ holds. 
$\square$
\end{theorem}

\begin{theorem}\label{semipositive}
Let $f : X \longrightarrow S$ be projective family such that 
$X$ and $S$ are smooth. 
Let $S^{\circ}$ be a nonempty Zariski open subset such that 
$f$ is smooth over $S^{\circ}$.
Then $K_{X/S}$ has a relative AZD $h$ over $S^{\circ}$ such that 
$\Theta_{h}$ is semipositive on $X$. 

And $F_{m} := f_{*}{\cal O}_{X}(mK_{X/S}))$ carries a continuous hermitian metric 
$h_{F_{m}}$ with Nakano semipositive curvature in the sense of current over
$S^{\circ}$. 

Let $x\in S- S^{\circ}$ be a point and let $\sigma$ be a local holomorphic
section of $F_{m}$ on a neighbourhood $U$ of $x$. 
Then $\sqrt{-1}\bar{\partial}\partial\log h_{F_{m}}(\sigma ,\sigma )$ extends as a closed positive current across 
$(S - S^{\circ})\cap U$.
$\square$
\end{theorem}

\begin{corollary}\label{invariance}
Let $f : X \longrightarrow S$ be a smooth projective family. 
Then $P_{m}(X_{s}) = \dim H^{0}(X_{s},{\cal O}_{X_{s}}(mK_{X_{s}}))$
is independent of $s\in S$. $\square$
\end{corollary}

After I completed writing this work, I have received a preprint of 
Bo Berndtsson \cite{b3}, which proved Theorem \ref{v}, under the assumption that 
$h$ is $C^{\infty}$. 
His proof is more computational than the one here and works also for  
smooth proper K\"{a}hler morphisms. 
But it looks quite different from 
his proof of Theorem \ref{b} which is very  ingeneous and beautiful. 
Also it is not clear whether his proof works also for singular $h$
although it seems not to be difficult at least for projective morphisms. 

The proof presented here is very simple and based on the beautiful  
proof of Theorem \ref{b} in \cite{b}. 
I would like to thank Professor Ohsawa for stimulating discussion.  

\section{Preliminaries}
\subsection{Singular hermitian metrics}\label{singh}
In this subsection $L$ will denote a holomorphic line bundle on a complex manifold $M$. 
\begin{definition}\label{singhm}
A  singular hermitian metric $h$ on $L$ is given by
\[
h = e^{-\varphi}\cdot h_{0},
\]
where $h_{0}$ is a $C^{\infty}$-hermitian metric on $L$ and 
$\varphi\in L^{1}_{loc}(M)$ is an arbitrary function on $M$.
We call $\varphi$ a  weight function of $h$. $\square$ 
\end{definition}
The curvature current $\Theta_{h}$ of the singular hermitian line
bundle $(L,h)$ is defined by
\[
\Theta_{h} := \Theta_{h_{0}} + \sqrt{-1}\partial\bar{\partial}\varphi ,
\]
where $\partial\bar{\partial}$ is taken in the sense of a current.
The $L^{2}$-sheaf ${\cal L}^{2}(L,h)$ of the singular hermitian
line bundle $(L,h)$ is defined by
\[
{\cal L}^{2}(L,h) := \{ \sigma\in\Gamma (U,{\cal O}_{M}(L))\mid 
\, h(\sigma ,\sigma )\in L^{1}_{loc}(U)\} ,
\]
where $U$ runs over the  open subsets of $M$.
In this case there exists an ideal sheaf ${\cal I}(h)$ such that
\[
{\cal L}^{2}(L,h) = {\cal O}_{M}(L)\otimes {\cal I}(h)
\]
holds.  We call ${\cal I}(h)$ the {\bf multiplier ideal sheaf} of $(L,h)$.
If we write $h$ as 
\[
h = e^{-\varphi}\cdot h_{0},
\]
where $h_{0}$ is a $C^{\infty}$ hermitian metric on $L$ and 
$\varphi\in L^{1}_{loc}(M)$ is the weight function, we see that
\[
{\cal I}(h) = {\cal L}^{2}({\cal O}_{M},e^{-\varphi})
\]
holds.
For $\varphi\in L^{1}_{loc}(M)$ we define the multiplier ideal sheaf of $\varphi$ by 
\[
{\cal I}(\varphi ) := {\cal L}^{2}({\cal O}_{M},e^{-\varphi}).
\] 
\begin{example}
Let $\sigma\in \Gamma (X,{\cal O}_{X}(L))$ be the global section. 
Then 
\[
h := \frac{1}{\mid\sigma\mid^{2}} = \frac{h_{0}}{h_{0}(\sigma ,\sigma)}
\]
is a singular hemitian metric on $L$, 
where $h_{0}$ is an arbitrary $C^{\infty}$-hermitian metric on $L$
(the right hand side is ovbiously independent of $h_{0}$).
The curvature $\Theta_{h}$ is given by
\[
\Theta_{h} = 2\pi\sqrt{-1}(\sigma )
\]
where $(\sigma )$ denotes the current of integration over the 
divisor of $\sigma$. $\square$ 
\end{example}
\begin{definition}\label{pe}
$L$ is said to be {\bf pseudoeffective}, if there exists 
a singular hermitian metric $h$ on $L$ such that 
the curvature current 
$\Theta_{h}$ is a closed positive current.
Also a singular hermitian line bundle $(L,h)$ is said to be {\bf pseudoeffective}, 
if the curvature current $\Theta_{h}$ is a closed positive current. $\square$
\end{definition}
\subsection{Analytic Zariski decompositions}
In this subsection we shall introduce the notion of analytic Zariski decompositions. 
By using analytic Zariski decompositions, we can handle  big line bundles
like  nef and big line bundles.
\begin{definition}\label{defAZD}
Let $M$ be a compact complex manifold and let $L$ be a holomorphic line bundle
on $M$.  A singular hermitian metric $h$ on $L$ is said to be 
an analytic Zariski decomposition, if the followings hold.
\begin{enumerate}
\item $\Theta_{h}$ is a closed positive current,
\item for every $m\geq 0$, the natural inclusion
\[
H^{0}(M,{\cal O}_{M}(mL)\otimes{\cal I}(h^{m}))\rightarrow
H^{0}(M,{\cal O}_{M}(mL))
\]
is an isomorphim. $\square$
\end{enumerate}
\end{definition}
\begin{remark} If an AZD exists on a line bundle $L$ on a smooth projective
variety $M$, $L$ is pseudoeffective by the condition 1 above. $\square$
\end{remark}

\begin{theorem}(\cite{tu,tu2})
 Let $L$ be a big line  bundle on a smooth projective variety
$M$.  Then $L$ has an AZD. 
\end{theorem}
As for the existence for general pseudoeffective line bundles, 
now we have the following theorem.
\begin{theorem}(\cite[Theorem 1.5]{d-p-s})\label{AZD}
Let $X$ be a smooth projective variety and let $L$ be a pseudoeffective 
line bundle on $X$.  Then $L$ has an AZD.
\end{theorem}
{\bf Proof of Theorem \ref{AZD}}.
Although the proof is in \cite{d-p-s}, 
we shall give a proof here, because we shall use it afterwards. 

 Let  $h_{0}$ be a fixed $C^{\infty}$-hermitian metric on $L$.
Let $E$ be the set of singular hermitian metric on $L$ defined by
\[
E = \{ h ; h : \mbox{lowersemicontinuous singular hermitian metric on $L$}, 
\]
\[
\hspace{70mm}\Theta_{h}\,
\mbox{is positive}, \frac{h}{h_{0}}\geq 1 \}.
\]
Since $L$ is pseudoeffective, $E$ is nonempty.
We set 
\[
h_{L} = h_{0}\cdot\inf_{h\in E}\frac{h}{h_{0}},
\]
where the infimum is taken pointwise. 
The supremum of a family of plurisubharmonic functions 
uniformly bounded from above is known to be again plurisubharmonic, 
if we modify the supremum on a set of measure $0$(i.e., if we take the uppersemicontinuous envelope) by the following theorem of P. Lelong.

\begin{theorem}(\cite[p.26, Theorem 5]{l})
Let $\{\varphi_{t}\}_{t\in T}$ be a family of plurisubharmonic functions  
on a domain $\Omega$ 
which is uniformly bounded from above on every compact subset of $\Omega$.
Then $\psi = \sup_{t\in T}\varphi_{t}$ has a minimum 
uppersemicontinuous majorant $\psi^{*}$  which is plurisubharmonic.
We call $\psi^{*}$ the uppersemicontinuous envelope of $\psi$. 
\end{theorem}
\begin{remark} In the above theorem the equality 
$\psi = \psi^{*}$ holds outside of a set of measure $0$(cf.\cite[p.29]{l}). 
\end{remark}

By Theorem 2.3,we see that $h_{L}$ is also a 
singular hermitian metric on $L$ with $\Theta_{h}\geq 0$.
Suppose that there exists a nontrivial section 
$\sigma\in \Gamma (X,{\cal O}_{X}(mL))$ for some $m$ (otherwise the 
second condition in Definition 2.3 is empty).
We note that  
\[
\log \mid\sigma\mid^{\frac{2}{m}}
\]
gives the weight of a singular hermitian metric on $L$ with curvature 
$2\pi m^{-1}(\sigma )$, where $(\sigma )$ is the current of integration
along the zero set of $\sigma$. 
By the construction we see that there exists a positive constant 
$c$ such that  
\[
\frac{h_{0}}{\mid\sigma\mid^{\frac{2}{m}}} \geq c\cdot h_{L}
\]
holds. 
Hence
\[
\sigma \in H^{0}(X,{\cal O}_{X}(mL)\otimes{\cal I}_{\infty}(h_{L}^{m}))
\]
holds.  
Hence in praticular
\[
\sigma \in H^{0}(X,{\cal O}_{X}(mL)\otimes{\cal I}(h_{L}^{m}))
\]
holds.  
 This means that $h_{L}$ is an AZD of $L$. $\square$ 
\vspace{10mm} 
\begin{remark}
By the above proof we have that for the AZD $h_{L}$ constructed 
as above
\[
H^{0}(X,{\cal O}_{X}(mL)\otimes{\cal I}_{\infty}(h_{L}^{m}))
\simeq 
H^{0}(X,{\cal O}_{X}(mL))
\]
holds for every $m$. $\square$
\end{remark}
It is easy to see that the multiplier ideal sheaves 
of $h_{L}^{m}(m\geq 1)$ constructed in the proof of
 Theorem 2.2 are independent of 
the choice of the $C^{\infty}$-hermitian metric $h_{0}$.
We call the AZD constructed as in the proof of Theorem \ref{AZD}  {\bf a canonical 
AZD} of $L$. 

\subsection{$L^{2}$-extension theorem}

\begin{theorem}(\cite[p.200, Theorem]{o-t})\label{o-t}
Let $X$ be a Stein manifold of dimension $n$, $\psi$ a plurisubharmonic 
function on $X$ and $s$ a holomorphic function on $X$ such that $ds\neq 0$ 
on every branch of $s^{-1}(0)$.
We put $Y:= s^{-1}(0)$ and 
$Y_{0} := \{ x\in Y; ds(x)\neq 0\}$.
Let $g$ be a holomorphic $(n-1)$-form on $Y_{0}$ 
with 
\[
c_{n-1}\int_{Y_{0}}e^{-\psi}g\wedge\bar{g} < \infty ,
\]
where $c_{k}= (-1)^{\frac{k(k-1)}{2}}(\sqrt{-1})^{k}$. 
Then there exists a holomorphic $n$-form $G$ on 
$X$ such that 
\[
G(x) = g(x)\wedge ds(x) 
\]
on $Y_{0}$ and 
\[
c_{n}\int_{X}e^{-\psi}(1+\mid s\mid^{2})^{-2}G\wedge\bar{G} 
\leq 1620\pi c_{n-1}\int_{Y_{0}}e^{-\psi}g\wedge\bar{g}. 
\] $\square$
\end{theorem}

For a extension from an arbitrary dimensional submanifold, 
T. Ohsawa extened Theorem \ref{o-t} in the following way. 

Let $M$ be a complex manifold of dimension $n$ and let $S$ be a closed complex submanifold of $M$. 
Then we consider a class of continuous function $\Psi : M\longrightarrow [-\infty , 0)$  such that  
\begin{enumerate}
\item $\Psi^{-1}(-\infty ) \supset S$,
\item if $S$ is $k$-dimensional around a point $x$, there exists a local 
coordinate $(z_{1},\ldots ,z_{n})$ on a neighbourhood of $x$ such that 
$z_{k+1} = \cdots = z_{n} = 0$ on $S\cap U$ and 
\[
\sup_{U\backslash S}\mid \Psi (z)-(n-k)\log\sum_{j=k+1}^{n}\mid z_{j}\mid^{2}\mid < \infty .
\]
\end{enumerate} 
The set of such functions $\Psi$ will be denoted by $\sharp (S)$. 

For each $\Psi \in \sharp (S)$, one can associate a positive measure 
$dV_{M}[\Psi ]$ on $S$ as the minimum element of the 
partially ordered set of positive measures $d\mu$ 
satisfying 
\[
\int_{S_{k}}f\, d\mu \geqq 
\overline{\lim}_{t\rightarrow\infty}\frac{2(n-k)}{v_{2n-2k-1}}
\int_{M}f\cdot e^{-\Psi}\cdot \chi_{R(\Psi ,t)}dV_{M}
\]
for any nonnegative continuous function $f$ with 
$\mbox{supp}\, f\subset\subset M$.
Here $S_{k}$ denotes the $k$-dimensional component of $S$,
$v_{m}$ denotes the volume of the unit sphere 
in $\mbox{\bf R}^{m+1}$ and 
$\chi_{R(\Psi ,t)}$ denotes the characteristic funciton of the set 
\[
R(\Psi ,t) = \{ x\in M\mid -t-1 < \Psi (x) < -t\} .
\]

Let $M$ be a complex manifold and let $(E,h_{E})$ be a holomorphic hermitian vector 
bundle over $M$. 
Given a positive measure $d\mu_{M}$ on $M$,
we shall denote $A^{2}(M,E,h_{E},d\mu_{M})$ the space of 
$L^{2}$ holomorphic sections of $E$ over $M$ with respect to $h_{E}$ and 
$d\mu_{M}$. 
Let $S$ be a closed  complex submanifold of $M$ and let $d\mu_{S}$ 
be a positive measure on $S$. 
The measured submanifold $(S,d\mu_{S})$ is said to be a set of 
interpolation for $(E,h_{E},d\mu_{M})$, or for the 
sapce $A^{2}(M,E,h_{E},d\mu_{M})$, if there exists a bounded linear operator
\[
I : A^{2}(S,E\mid_{S},h_{E},d\mu_{S})\longrightarrow A^{2}(M,E,h_{E},d\mu_{M})
\]
such that $I(f)\mid_{S} = f$ for any $f$. 
$I$ is called an interpolation operator.
The following theorem is crucial.

\begin{theorem}(\cite[Theorem 4]{o})\label{extension}
Let $M$ be a complex manifold with a continuous volume form $dV_{M}$,
let $E$ be a holomorphic vector bundle over $M$ with $C^{\infty}$-fiber 
metric $h_{E}$, let $S$ be a closed complex submanifold of $M$,
let $\Psi\in \sharp (S)$ and let $K_{M}$ be the canonical bundle of $M$.
Then $(S,dV_{M}(\Psi ))$ is a set of interpolation 
for $(E\otimes K_{M},h_{E}\otimes (dV_{M})^{-1},dV_{M})$, if 
the followings are satisfied.
\begin{enumerate}
\item There exists a closed set $X\subset M$ such that 
\begin{enumerate}
\item $X$ is locally negligble with respect to $L^{2}$-holomorphic functions, i.e., 
for any local coordinate neighbourhood $U\subset M$ and for any $L^{2}$-holomorphic function $f$ on $U\backslash X$, there exists a holomorphic function 
$\tilde{f}$ on $U$ such that $\tilde{f}\mid U\backslash X = f$.
\item $M\backslash X$ is a Stein manifold which intersects with every component of $S$. 
\end{enumerate}
\item $\Theta_{h_{E}}\geqq 0$ in the sense of Nakano,
\item $\Psi \in \sharp (S)\cap C^{\infty}(M\backslash S)$,
\item $e^{-(1+\epsilon )\Psi}\cdot h_{E}$ has semipositive 
curvature in the sense of Nakano for every $\epsilon \in [0,\delta]$ 
for some $\delta > 0$.
\end{enumerate}
Under these conditions, there exists a constant $C$ and an interpolation operator 
from $A^{2}(S,E\otimes K_{M}\mid_{S},h\otimes (dV_{M})^{-1}\mid_{S},dV_{M}[\Psi ])$
to $A^{2}(M,E\otimes K_{M},h\otimes (dV_{M})^{-1}.dV_{M})$ whose 
norm does not exceed $C\delta^{-3/2}$.
If $\Psi$ is plurisubharmonic, the interpolation operator can be chosen 
so that its norm is less than $2^{4}\pi^{1/2}$. $\square$
\end{theorem}

The above theorem can be generalized to the case that 
$(E,h_{E})$ is a singular hermitian line bundle with semipositive
curvature current  (we call such a singular hermitian line 
bundle $(E,h_{E})$ a {\bf pseudoeffective singular hermitian line bundle}) as was remarked in \cite{o}. 

\begin{lemma}\label{extension2}
Let $M,S,\Psi ,dV_{M}, dV_{M}[\Psi], (E,h_{E})$ be as in Theorem \ref{extension} 
Let $(L,h_{L})$ be a pseudoeffective singular hermitian line 
bundle on $M$. 
Then $S$ is a set of interpolation for 
$(K_{M}\otimes E\otimes L,dV_{M}^{-1}\otimes h_{E}\otimes h_{L})$.  $\square$
\end{lemma}

\section{Proof of Theorem \ref{v}}
Let $f : X \longrightarrow S$ be a projective family. 
Since the statement is local we may assume that $S$ is the unit open disk 
$\Delta$ in {\bf C}.
Let us consider the family 
\[
f : X(\varepsilon ) \longrightarrow \Delta
\]
be the family such that 
\[
X(\varepsilon )_{t} := f^{-1}(\mbox{the $\varepsilon$-neighbourhood of $t$ in $\Delta$})
\]
where $\varepsilon$ means the $\varepsilon$-neighbourhood with repsect to 
the Poincar\'{e} metric on $\Delta$. 
Then $X(\varepsilon )$ is a pseudoconvex domain in $X\times \Delta$. 
Then the Bergman kernel 
\[
K(X(\varepsilon )_{t},K_{X}+L,h_{L}), 
\]
satisfies that 
\[
\sqrt{-1}\partial\bar{\partial}\log K(X(\varepsilon )_{t},K_{X}+L,h_{L}), 
 \geqq 0
\]
on $X(\varepsilon )$ by Theorem \ref{b}.
Since 
\[
\lim_{\varepsilon\downarrow 0}(\pi\varepsilon^{2})\cdot K(X(\varepsilon )_{t},K_{X}+L,h_{L})
= K(X_{t},K_{X_{t}} + L, h\mid X_{t})
\]
holds.  

In fact, if we consider the family 
\[
\pi_{\varepsilon ,t} : X(\varepsilon ) \longrightarrow \Delta (t,\varepsilon ),
\]
where 
\[
\Delta (t,\varepsilon ) := \{ t^{\prime}\in\mbox{\bf C} \mid 
\mid t^{\prime} - t\mid < \varepsilon \}
\]
as a family over the unit open disk $\Delta$ in {\bf C} by 
\[
t^{\prime}\mapsto \varepsilon^{-1}(t^{\prime} - t),
\]
the limit as $\varepsilon \downarrow 0$ is nothing but the 
trivial family $X_{t}\times \Delta$. 

This completes the proof of Theorem \ref{v}. $\square$ 

\section{Proof of Theorem \ref{semipositive}}

\subsection{Dynamical construction of an AZD}\label{Dy}

Let $X$ be a smooth projective variety and let 
$K_{X}$ be the canonical line bundle of $X$.
Let $n$ denote the dimension of $X$.
We shall assume that $K_{X}$ is pseudoeffective.  
Then by Theorem , $K_{X}$ admits an AZD $h$. 

Let $A$ be a sufficiently ample line bundle on $X$ 
such that for every pseudoeffective singular hermitian 
line bundle $(L,h_{L})$
\[
{\cal O}_{X}(A+L)\otimes{\cal I}(h_{L})
\]
and 
\[
{\cal O}_{X}(K_{X}+A+L)\otimes{\cal I}(h_{L})
\]
are globally generated. 
This is possible by \cite[p. 667, Proposition 1]{si}. 

Let $h_{A}$ be a $C^{\infty}$ hermitian metric on $A$
 with strictly positive curvature.  
For $m\geq 0$, let $h_{m}$ be the singular hermitian metrics
on $A + mK_{X}$ constructed as follows. 
Let $h_{0}$ be a $C^{\infty}$-hermitian metric $h_{A}$ on 
$A$ with strictly positive curvature. 
Suppose that $h_{m-1} (m\geq 1)$ has been constructed. 
We set 
\begin{eqnarray*}
K_{m} & :=  &K(X,m(K_{X},h)+A, h_{m-1}) \\
& = & \sup \{ \mid \sigma\mid^{2} \mid 
(\sqrt{-1})^{n^{2}}\int_{X}h_{m-1}\sigma\wedge \bar{\sigma} \leqq 1, 
\sigma \in H^{0}(X,{\cal O}_{X}(mK_{X}+A)\otimes {\cal I}(h^{m}))\} .
\end{eqnarray*}
And we define the singular hermitian metric $h_{m}$ 
on $A + mK_{X}$ by  
\[
h_{m} := K_{m}^{-1}.
\]
It is clear that $K_{m}$ has semipositive curvature 
in the sense of currents. 
We note that for every $x\in X$
\[
K_{m}(x) = \sup \{ \mid\sigma\mid^{2}(x) ;  
\sigma\in \Gamma (X,{\cal O}_{X}(A+mK_{X}\otimes{\cal I}(h^{m}))), 
\int_{X}h_{m-1}\cdot \mid\sigma\mid^{2} = 1\}
\]
holds by definition (cf. \cite[p.46, Proposition 1.4.16]{kr}).

We set 
\[
\nu := \overline{\lim}_{m\rightarrow\infty}\frac{\log \dim H^{0}(X,{\cal O}_{X}(A + mK_{X})\otimes{\cal I}(h^{m}))}{\log m}
\]
and call it the numerical Kodaira dimension of $(K_{X},h)$. 

\begin{proposition}\label{dynamical}(cf. \cite{tu3})
\[
K_{\infty}:= \overline{\lim}_{m\rightarrow\infty}\sqrt[m]{(m!)^{-\nu}K_{m}}
\]
exists and 
\[
h_{\infty}:= 1/K_{\infty}
\]
is an AZD of $K_{X}$. $\square$ 
\end{proposition}
{\bf Proof of Proposition \ref{dynamical}}.
There exists a positive constant $C$ such that 
\[
h^{0}(X,{\cal O}_{X}(mK_{X} + A)\otimes{\cal I}(h)) \leqq C m^{\nu}
\]
holds. 
Let $dV$ be a fixed $C^{\infty}$ volume form on $X$. 
Then by the submeanvalue inequality of plurisubharmonic function, we 
see that by induction there exists a positive constant $C_{1}$ such that 
\begin{equation}\label{upper}
h_{A}\cdot K_{m} \leqq C_{1}^{m}\cdot (m!)^{\nu}dV^{m}
\end{equation}
holds.

Let $h$ be an AZD of $K_{X}$. 

\[
K_{m}(x) = \sup \{ \mid\sigma\mid^{2}(x) ;  
\sigma\in \Gamma (X,{\cal O}_{X}(A+mK_{X})), 
\int_{X}h_{m-1}\cdot \mid\sigma\mid^{2} = 1\}.
\]

Let $x \in X$ be a point. 
Let $H$ be a sufficiently ample divisor on $X$. 
Let $H_{1}, \cdots ,H_{n-\nu}$ be a general member of $\mid H\mid$ 
containing $x$.
And let $V := H_{1}\cap \cdots \cap H_{n-\nu}$.
Then the restriction morphism 
\[
H^{0}(X,{\cal O}_{X}(mK_{X}+A )\otimes{\cal I}(h^{m})) \rightarrow 
H^{0}(V,{\cal O}_{V}(mK_{X}+A)\otimes{\cal I}(h^{m}))
\]
is injective. 
Hence $(K_{X},h)\mid V$ is big.  
By Kodaira's lemma, 
$h\mid V$ is dominated by a singular hermitian metric $h^{\prime}$ 
of $L\mid V$ such that 
$\Theta_{h^{\prime}}$ is strictly positive on $V$. 
For $0 < \varepsilon < 1$ we set 
\[
h_{V,\varepsilon}:= (h\mid V)^{1-\varepsilon}\cdot (h^{\prime})^{\varepsilon}.
\]
Then $h\mid V < h_{V,\varepsilon}$ holds.

Suppose that 
\begin{equation}\label{ind}
h_{A}\cdot K_{m} \geqq C(m)(m!)^{\nu}\cdot h_{V,\varepsilon}^{-m}
\end{equation}
holds on $V$ for some positive constant $C(m)$.

Then by the $L^{2}$-extension theorem, there exists a positive constant $C$
such that 
\begin{equation}\label{ind2}
h_{A}\cdot K_{m+1} \geqq C\cdot C(m)(m!)^{\nu}\cdot (m+1)^{\nu}\cdot h_{V,\varepsilon}^{-(m+1)}
\end{equation} 
holds on $V$, where $C$ is a positive constant independent of $m$. 

Here we have applied the $L^{2}$-extension theorem, first to 
the extension from a point $x \in V$ and the second to the extension 
from $V$ to $X$. 
The constant $(m+1)^{\nu}$ appears simply because $h_{V,\varepsilon}$ is 
strictly positive, hence we can take local frame $\mbox{\bf e}$ of $K_{X}$ 
around  $x\in V$ and coordinate $z_{1},\cdots ,z_{\nu}$ so that 
\[
h_{V,\varepsilon}(\mbox{\bf e},\mbox{\bf e})
= (1 - \parallel z\parallel^{2})h (\mbox{\bf e},\mbox{\bf e})(x) + o(\parallel z\parallel^{2})
\]
holds (cf\cite[p.105,(1,11)]{ti}). 

Hence combining (\ref{ind}) and (\ref{ind2}) we have that 
there exists a positive constant $C_{2}$ indendent of $m$ such that 
\begin{equation}\label{lower}
h_{A}\cdot K_{m} \geqq C_{2}^{m}(m!)^{-\nu}\cdot h_{V,\varepsilon}^{-m}
\end{equation}
holds on $V$. 

By (\ref{upper}) and (\ref{lower}), moving $x$ and $V$ and 
letting $\varepsilon$ tend to $0$, we have that
\[
K_{\infty}:= \overline{\lim}_{m\rightarrow\infty}\sqrt[m]{(m!)^{-\nu}K_{m}}
\]
exists and 
\[
h_{\infty}:= 1/K_{\infty}
\]
is an AZD of $K_{X}$. 

This completes the proof of Proposition \ref{dynamical}. $\square$ 

\subsection{Dynamical construction of AZD as a family}

Let $f : X \longrightarrow \Delta$ be a smooth projective family such that
$K_{X}$ is pseudoeffective.  
Let $n$ be the relative dimension of $f : X \longrightarrow \Delta$. 
Let $A$ be a sufficiently ample line bundle as in the last subsection
and let $h_{A}$ be a hermitian metric with strictly positive curvature
on $X$. 

Let $h_{0}$ be an AZD of $X_{0} = f^{-1}(0)$.
We define the sequence of Bergman kernels
\[
K_{m,0}:= K(X_{0},m(K_{X_{0}},h_{0})+A,h_{m-1,0})
\]
starting
\[
K_{0,0} = (h_{A}\mid X_{0})^{-1}
\]
as in the last subsection. 
Let $\{ \sigma^{(m)}_{0},\cdots ,\sigma^{(m)_{N(m)}}\}$ 
be a set of orthonormal basis of 
\[
H^{0}(X,{\cal O}_{X}(mK_{X_{0}}+A\mid X_{0})\otimes {\cal I}(h_{0}^{m}))
\]
with respect to the innner product
\[
(\sigma ,\sigma^{\prime}) := (\sqrt{-1})^{n^{2}}\int_{X_{0}}h_{m-1}\cdot\sigma\wedge\bar{\sigma}^{\prime}.
\]
Inductively we extend each $\sigma^{m}_{i}$ to 
\[
\tilde{\sigma}^{(m)}_{i} \in H^{0}(X,{\cal O}_{X}(mK_{X} + A)
\otimes {\cal I}(\tilde{h}_{m-1}))
\]
with the estimate
\[
\parallel\tilde{\sigma}^{(m)}_{i}\parallel^{2} 
= (\sqrt{-1})^{n^{2}}\int_{X}\tilde{h}_{m-1}\tilde{\sigma}^{(m)}_{i}\wedge\tilde{\sigma}^{(m)}_{i} \leqq C, 
\]
where $C$ is a positive constant indepdent of $m$ and $i$ and 
$\tilde{h}_{m-1}$ is defined inductively by 
\[
\tilde{h}_{0} = h_{A}
\]
and 
\[
\tilde{h}_{j} := \frac{1}{K_{j}}
\]
where
\[
\tilde{K}_{j}:=\sum_{i=0}^{N(j)}\mid\tilde{\sigma_{i}}^{(j)}\mid^{2}. 
\]
Then by the same argument as in Section \ref{Dy}, we see that 
\[
K_{\infty}:= \overline{\lim}_{m\rightarrow\infty}\sqrt[m]{(m!)^{-\nu}K_{m,0}}
\]
and 
\[
\tilde{K}_{\infty}:= \overline{\lim}_{m\rightarrow\infty}\sqrt[m]{(m!)^{-\nu}\tilde{K}_{m}}
\]
exists and both nonzero.  
Hence we see that an AZD of $K_{X}$ restricts to an AZD of 
$K_{X_{0}}$. 

\begin{proposition}\label{family}
Let $f : X \longrightarrow \Delta$ be a smooth projective family such that 
$K_{X}$ is pseudoeffective. 
Then for an AZD $h$ of $K_{X}$.  The restriction $h\mid X_{t}$ is an 
AZD of $K_{X_{t}}$. $\square$ 
\end{proposition}
\begin{remark}
As in \cite{tu3}, by using Theorem \ref{extension}, 
Proposition \ref{family} implies the invariance of plurigenra,
Corollary \ref{invariance}.  $\square$
\end{remark}

Now we shall prove Theorem \ref{semipositive}. 

Let $f : X \longrightarrow \Delta$ be a smooth projective family
over the unit open disk $\Delta$ with center $O$.
Let $A$ be a sufficiently ample line bundle on $X$ and 
let $h_{A}$ be a $C^{\infty}$ hermitian metric on $A$.  
Let $h$ be an AZD of $K_{X}$ constructed as in Theorem \ref{AZD}. 
Then by Proposition \ref{family}, we have that the restriction of 
$h$ to $X_{t}:= f^{-1}(t)$ is an AZD of $K_{X_{t}}$ for every 
$t\in\Delta$. 
Let $\nu_{t}$ be the numerical Kodaira dimension of $(K_{X_{t}},h\mid X_{t})$.
By Proposition \ref{family}, we see that $\nu_{t}$ is independent of 
$t\in \Delta$.  Hence we shall denote $\nu_{t}$ simply by $\nu$. 

Let us perform the dynamical construction of AZD as in Section \ref{Dy}.
Namely for every $t\in \Delta$, we start from $h_{A}\mid X_{t}$, by 
induction we define the hermitian metric $h_{m,t}$ as in 
Section \ref{Dy}.  
Then 
\[
h_{\infty ,t} := \liminf_{m\rightarrow\infty}\sqrt[m]{(m!)^{\nu}h_{m,t}}
\]
is an AZD of $K_{X_{t}}$. 
By Theorem \ref{v}, we see that the singular hermitian metric
$h_{m}$ on $mK_{X/\Delta}+A$ defined by 
\[
h_{m}\mid X_{t} = h_{m,t}  (t\in \Delta )
\]
has semipositive curvature in the sense of current on $X$. 
Then by the construction the singular hermitian metric $h_{\infty}$ 
on $K_{X/\Delta}$ defined by 
\[
h_{\infty}\mid X_{t} = h_{\infty ,t}
\]
has semipositve curvature in the sense of current on $X$. 
Hence by the construction $h_{\infty}$ is an AZD of $K_{X/\Delta}$. 
This completes the proof of Theorem \ref{semipositive} except the last 
assertion. 

The last assertion follows from the fact that the neighbourhood of 
the singular fiber is a manifold.  Hence the proof of Theorem \ref{v} 
implies the assertion. 
 $\square$. 

\subsection{Case of general type}

Let $m$ be a positive integer.
For a section $\eta \in H^{0}(X,{\cal O}_{X}(mK_{X}))$ we define  
a nonnegative number $\parallel\eta\parallel_{\frac{1}{m}}$ by 
\[
\parallel\eta\parallel_{\frac{1}{m}} = \mid \int_{X}(\eta\wedge\bar{\eta})^{\frac{1}{m}}\mid^{\frac{m}{2}}.
\]
Then $\eta\mapsto \parallel\eta\parallel_{\frac{1}{m}}$ is a continuous pseudonorm 
on $H^{0}(X,{\cal O}_{X}(mK_{X}))$, i.e.,  it is a continuous and has the properties :
\begin{enumerate}
\item $\parallel\eta\parallel_{\frac{1}{m}} = 0 \Leftrightarrow \eta = 0$,
\item $\parallel\lambda \eta\parallel_{\frac{1}{m}} = \mid\lambda\mid\cdot
\parallel\eta\parallel_{\frac{1}{m}} $ holds for all $\lambda\in \mathbb{C}$. 
\end{enumerate}
But it is not a norm on $H^{0}(X,{\cal O}_{X}(mK_{X}))$ except $m=1$.
We define a continuous section $K_{m}$ of 
$(K_{X}\otimes \bar{K}_{X})^{\otimes m}$ 
\[
K^{NS}_{m}(x):= \sup\{ \mid \eta (x)\mid^{2}\,\,\,; 
\,\,\,\parallel \eta \parallel_{\frac{1}{m}} = 1\}  \,\,\,\,\,\,\, (x\in X),  
\]
where $\mid\eta (x)\mid^{2} = \eta (x)\otimes \bar{\eta}(x)$.
We call $K^{NS}_{m}$ {\bf the $m$-th Narashimhan-Simha potential} of $X$.
We define the singular hermitian metric $h^{NS}_{m}$ on 
$H^{0}(X,{\cal O}_{X}(mK_{X}))$ by 
\[
h^{NS}_{m} := \frac{1}{K^{NS}_{m}}.
\]
We call $h^{NS}_{m}$ {\bf the Narashimhan-Simha metric} on $mK_{X}$. 
This metric is introduced by Narashimhan and Simha for smooth
 canonically polarized variety to study the moduli space of canonically polarized variety (\cite{n-s}).  
We note that the singularities of $h_{m}$ is located exactly on the support of 
the base locus of $\mid mK_{X}\mid$. 
Hence in the case of canonically polarized variety, $h_{m}$ is a nonsingular continuous metric on $mK_{X}$ for every sufficiently large $m$.  
Since $K^{NS}_{m}$ is locally the supremum of the  family of powers of absolute value of  holomorphic functions we see that the curvature
\[
\Theta_{h^{NS}_{m}} := \sqrt{-1}\partial\bar{\partial}\log K^{NS}_{m}
\] 
is a closed positive current. 

\begin{theorem}(\cite[Main Theorem]{tu4})\label{NS}
Let $f : X \longrightarrow S$ be a flat projective family of varieties with only canonical singularities over a complex manifold $S$. 
Let $h^{NS}_{m}$ be the Narashimhan-Simha singular hermitian metric  on $mK_{X/S}$.
 
Then $h^{NS}_{m}$ has semipositive curvature in the sense of current on $X$. $\square$  
\end{theorem}

The dynamical construction of AZD in Section \ref{Dy} works without 
an ample line bundle $A$ when the manifold is of general type and 
we can make the construction canonical. 

The construction is as follows. 
Let $m_{0}$ be a sufficiently large  positive integer such that $\mid m_{0}K_{X}\mid$ gives a birational embedding of $X$. 
Let $h_{m_{0}}^{NS}$ be the Narashimhan-Simha metric on  $m_{0}K_{X}$. 
Then starting $h_{m_{0}}^{NS}$ we may construct an AZD as in Section \ref{Dy}. 
In fact one may easily seen that for a sufficiently large $m_{0}$, the 
movable part of $\mid m_{0}K_{X}\mid$ dominates an ample divisor $A$ 
used in Section \ref{Dy}. 

Let $f : X \longrightarrow S$ be a smooth projective family of 
manifolds of general type. 

Then by Theorem \ref{NS}, we may construct an AZD of 
$K_{X/S}$ which is canonical.

Hence we obtain the following theorem. 
\begin{theorem}\label{general}
Let $f : X \longrightarrow S$ be a smooth projective family of 
manifolds of general type. 
Then for every positive integer $m$, there exists a functorial $C^{0}$ hermitian metric $h_{F_{m}}$ on the vector bundle $F_{m}:= f_{*}{\cal O}_{X/S}(mK_{X/S})$ with semipositive curvature current in the sense of Nakano. 

Here ``functorial'' means that $h_{F_{m}}$ only depends on 
the birational moduli map ;  
$\mbox{\bf birmod} : S \longrightarrow {\cal M}_{bir}$, 
where ${\cal M}_{bir}$ denote the (set theoretic) moduli space 
of the birational equivalence classes of projective varieties of general type
(after fixing $m_{0}$).  $\square$
\end{theorem}

As in \cite{tu4}, using Theorem \ref{general}, we may prove  
the quasiprojectivity of the moduli space of canonically polarized 
varieties.

Author's address\\
Hajime Tsuji\\
Department of Mathematics\\
Sophia University\\
7-1 Kioicho, Chiyoda-ku 102-8554\\
Japan \\
e-mail address: tsuji@mm.sophia.ac.jp

\end{document}